\documentclass[11pt]{article}
\usepackage{graphicx}
\usepackage{amsmath}
\usepackage{amsfonts}
\usepackage{amssymb}

\def\d{\mathrm d}
\def\e{\mathrm e}

\def\limfdd{\renewcommand{\arraystretch}{0.5}
\begin{array}[t]{c}
\stackrel{\rm fdd}{\longrightarrow} \\
\end{array}\renewcommand{\arraystretch}{1}}

\def\i{\mathrm i}

\newtheorem{thm}{Theorem}[section]

\newtheorem{prop}[thm]{Proposition}

\newtheorem{defn}[thm]{Definition}

\newtheorem{rem}[thm]{Remark}

\numberwithin{equation}{section}

\newcommand{\R}{\mathbb R}

\newcommand{\E}{\mathbb E}
\newcommand{\C}{\mathbb C}

\newcommand{\N}{\mathbb N}

\newcommand{\Z}{\mathbb Z}

\newcommand{\1}{\mathbf 1}

\def\eqfdd{\renewcommand{\arraystretch}{0.5}
\begin{array}[t]{c}
\stackrel{\rm fdd}{=} \\
\end{array}\renewcommand{\arraystretch}{1}}

\def\eqd{\renewcommand{\arraystretch}{0.5}
\begin{array}[t]{c}
\stackrel{\rm d}{=} \\
\end{array}\renewcommand{\arraystretch}{1}}

\newcommand{\limp}{\stackrel{p}{\longrightarrow}}

\newcommand{\DD}{\stackrel{D[0,1]}{\longrightarrow}}

\newcommand{\nn}{\nonumber}
\newcommand{\noi}{\noindent}

\def\limd{\renewcommand{\arraystretch}{0.5}
\begin{array}[t]{c}
\stackrel{\rm d}{\longrightarrow} \\
\end{array}\renewcommand{\arraystretch}{1}}

\def\limfdd{\renewcommand{\arraystretch}{0.5}
\begin{array}[t]{c}
\stackrel{\rm fdd}{\longrightarrow} \\
\end{array}\renewcommand{\arraystretch}{1}}

\def\limfdd{\renewcommand{\arraystretch}{0.5}
\begin{array}[t]{c}
\stackrel{\rm fdd}{\longrightarrow} \\
\end{array}\renewcommand{\arraystretch}{1}}

\begin{document}

\title
{Invariance Principles for Tempered Fractionally \\
Integrated Processes}

\author{Farzad Sabzikar$^{1}$  \  and \  Donatas Surgailis$^2$}
\date{\today \\  \small
\vskip.2cm
$^1$Iowa State University \ and \
$^2$Vilnius University}

\maketitle

\begin{abstract}

We discuss invariance principles 
for autoregressive tempered fractionally integrated moving averages 
in $\alpha$-stable $(1< \alpha \le 2)$ i.i.d. innovations and related tempered linear processes
with vanishing tempering parameter $\lambda \sim \lambda_*/N$. 
We show that the limit of the partial sums process takes a different form
in the weakly tempered ($\lambda_* = 0$), strongly tempered  ($\lambda_* = \infty$),
and moderately tempered ($0<\lambda_* < \infty$) cases. These results are used to  derive
the limit distribution of the OLS estimate of  AR(1) unit root with  weakly, strongly, and moderately tempered
moving average errors.


\end{abstract}

\smallskip
{\small

\noi {\it Keywords:} invariance principle; tempered linear process; autoregressive fractionally integrated moving average;
tempered fractional stable/Brownian motion;
tempered fractional unit root distribution;

}

\section{Introduction}

The present paper discusses  partial sums limits and invariance principles for tempered 
moving averages 
\begin{equation}\label{TL}
X_{d,\lambda}(t) =
\sum_{k=0}^{\infty} \e^{-\lambda k} b_d(k)\zeta (t-k), \qquad t \in \Z
\end{equation}
in i.i.d. innovation process $\{\zeta (t)\}$ with coefficients $b_d(k)$ regularly varying at infinity as $k^{d-1}$, viz.
\begin{equation} \label{bdk}
b_d (k)\ \sim \ \frac{c_d}{\Gamma(d)} \, k^{d-1}, \qquad k \to \infty,  \qquad c_d \ne 0,  \quad d \ne 0
\end{equation}
where $ d \in \R  $ is a real number, $d \ne -1,-2, \dots $ and
$\lambda >0$ is tempering parameter.  In addition to \eqref{bdk} we assume that
\begin{eqnarray}
&&\sum_{k=0}^\infty k^j b_d(k) = 0,  \quad 0 \le j  \le [-d],  \hskip.5cm
  -\infty  < d < 0,  \label{bdkneg} \\
&&\sum_{k=0}^\infty |b_d(k)| < \infty,  \quad  \sum_{k=0}^\infty b_d(k) \ne 0,  \qquad  d = 0 \label{bdk0}
\end{eqnarray}
An important example of such processes is
the two-parametric  class
ARTFIMA$(0,d,\lambda,0)$ of
tempered fractionally integrated processes, generalizing the well-known   ARFIMA$(0,d,0)$ class,  written as
\begin{equation}\label{ART}
X_{d,\lambda}(t) = (1- \e^{-\lambda} B)^{-d} \zeta(t)  =
\sum_{k=0}^{\infty} \e^{-\lambda k} \omega_{-d}(k)\zeta (t-k), \qquad t \in \Z
\end{equation}
with coefficients given by power expansion $(1- \e^{-\lambda} z)^{-d} = \sum_{k=0}^\infty \e^{-\lambda k} \omega_{-d}(k) z^k, |z| < 1$,  where $B x(t) = x(t-1)$ is the backward shift. Due to the presence of the exponential tempering factor $\e^{-\lambda k}$ the series
in \eqref{TL} and
\eqref{ART} absolutely converges a.s. under general assumptions on the innovations, and defines
a strictly stationary process. On the other hand, for $\lambda = 0$  
the corresponding stationary 
processes in \eqref{TL} and
\eqref{ART} exist under additional  conditions on the parameter $d$. See  Granger and Joyeux~\cite{Granger}, Hosking~\cite{Hosking},
Brockwell and Davis~\cite{BrockwellDavisTSTM}, Kokoszka and Taqqu~\cite{kokoszka}.
We also note (see e.g. \cite{koul}, Ch.~3.2) 
that the (untempered) linear process $X_{d,0}$ of \eqref{TL} 
with coefficients satisfying \eqref{bdk} for $0< d < 1/2 $ is said long memory,
while \eqref{bdk} and \eqref{bdkneg} for
$-1/2 < d < 0$ is termed negative memory and
\eqref{bdk0} short memory, respectively, parameter $d$ usually referred to as memory parameter.

The model in \eqref{ART}  appeared in  Giraitis et al.~\cite{giraitis}, which noted
that for small $\lambda >0$,  $X_{d,\lambda}$ has a covariance function
which resembles the covariance function of a long memory model for arbitrary large number of lags but
eventually decays exponentially fast.  \cite{giraitis} termed such behavior `semi long-memory'
and noted that it may have empirical  relevance for modelling of financial returns.
Giraitis et al.~\cite{giraitis2003a} propose
the semi-long memory ARCH($\infty$) model as a contiguous alternative to (pure) hyperbolic and
exponential decay which are often very hard to distinguish between in a finite sample. On the other side,
Meerschaert et al.~\cite{Meerschaertsabzikarkumarzeleki} effectively apply ARTFIMA$(0,d,\lambda,0)$ in
\eqref{ART} for modeling of turbulence in the Great Lakes region.

The present paper obtains limiting behavior of tempered linear processes in \eqref{TL} 
with small tempering parameter $\lambda = \lambda_N \to 0$ tending to zero together with
the sample size. The important statistic is the partial sums process
\begin{equation}\label{Spart}
S^{d,\lambda}_{N}(t) :=   \sum_{k=1}^{[Nt]}X_{d,\lambda}(k),  \qquad t\in [0,1]
\end{equation}
of $X_{d,\lambda}$ in \eqref{TL} with i.i.d. innovations $\{\zeta (t)\}$ in the domain of attraction
of $\alpha$-stable law, $1 < \alpha \le 2$.  Functional limit theorems for the partial sums
process play a crucial role in the R/S analysis, unit root testing,  change-point analysis and many other time series inferences.
See  Lo~\cite{Lo}, Phillips~\cite{Phillips}, Giraitis et al.~\cite{giraitis2003a},  Lavancier et al.~\cite{Lavancier} and the references
therein.

We prove that the limit behavior of \eqref{Spart} essentially depends on how fast $\lambda = \lambda_N $ tends to 0.
Assume that there exists the limit
\begin{equation}\label{lambdalim}
\lim_{N \to \infty}  N \lambda_N = \lambda_*  \in [0, \infty].
\end{equation}
Depending on the value of $\lambda_*$, the process  $X_{d,\lambda_N}$ will be called
{\it strongly tempered} if $\lambda_* =  \infty$,
{\it weakly tempered} if $\lambda_* =  0$, and
{\it moderately tempered} if $0 < \lambda_* < \infty $.  While the behavior of $S^{d,\lambda_N}_{N}$
in the strongly and weakly tempered cases is typical for short memory and long memory processes, respectively,
the moderately tempered decay $\lambda_N \sim \lambda_*/N, \lambda_* \in (0,\infty)$
leads to {\it tempered fractional stable motion of  second kind} (TFSM II) $Z^{I\!I}_{H,\alpha,\lambda_*}, H=  d + 1/\alpha >0$ 
defined as a stochastic integral
\begin{eqnarray}\label{T2}
Z^{I\!I}_{H,\alpha,\lambda}(t)&:=&\int_{\R} h_{H,\alpha,\lambda}(t;y) M_\alpha (\d y), \qquad t \in \R
\end{eqnarray}
with respect to $\alpha$-stable L\'evy process $M_\alpha $ with integrand
\begin{eqnarray}\label{hdef0}
h_{H,\alpha,\lambda}(t;y)&:=&(t-y)_+^{H - \frac{1}{\alpha}} \e^{-\lambda (t-y)_+} - (-y)_+^{H - \frac{1}{\alpha}} \e^{-\lambda (-y)_+} \\
&+&\lambda \int_{0}^{t} (s-y)_{+}^{H-\frac{1}{\alpha}} \e^{-\lambda(s-y)_{+}}\ \d s, \qquad y \in \R. \nn
\end{eqnarray}
TFSM II and its Gaussian counterpart {\it tempered fractional Brownian motion of  second kind} (TFBM II) 
were recently introduced in Sabzikar and Surgailis~\cite{SurgailisFarzadTFSMII},
the above processes being closely related to the {\it tempered fractional stable motion} (TFSM) and
the {\it tempered fractional Brownian motion} (TFBM)
defined in
Meerschaert and Sabzikar~\cite{Meerschaertsabzikar2} and
Meerschaert and Sabzikar~\cite{Meerschaertsabzikar}, respectively.
As shown in \cite{SurgailisFarzadTFSMII}, 
TFSM and TFSM II are different processes, 
especially striking are their differences  as $t \to \infty $.

As an application of our invariance principles we obtain the limit distribution
of the OLS estimator $\widehat \beta_N $ of the slope parameter in AR(1) model  
with tempered ARTFIMA$(0,d,\lambda_N,  0)$
errors and small tempering
parameter $\lambda_N\to 0$ satisfying \eqref{lambdalim},
under the null (unit root)
hypothesis $\beta = 1 $.
In the case of (untempered) ARFIMA$(0,d,0)$
error process with finite variance and  standardized i.i.d. innovations,
Sowell~\cite{sowell2} proved that the distribution of the normalized statistic
$N^{1 \wedge (1+2d)} (\widehat \beta_N - 1) $ tends to the so-called 
{\it fractional unit root distribution} written in terms of fractional Brownian motion
with parameter $H = d + \frac{1}{2}$. 
Sowell's~\cite{sowell2}  result extends the classical unit root distribution
for weakly dependent errors in Phillips~\cite{Phillips}
to  fractionally integrated error process, yielding drastically different limits
for $0< d < 1/2, d = 0$ and $-1/2 < d < 0$.

It turns out
that in the case of ARTFIMA$(0,d,\lambda_N,  0)$ error process
with $\lambda_N \sim \lambda_*/N $, the limit distribution
of $\widehat \beta_N$ depends on $\lambda_* \in [0,\infty]$ and $d$. Roughly speaking
(see Theorem~\ref{thm:ADDIDTIVE ARTFIMA TO REGRESSION MODEL1}  for precise formulation), in the moderately tempered case
$0<\lambda_* < \infty $ the  limit distribution of  $\widehat \beta_N$ 
writes similarly to  Sowell~\cite{sowell2} with FBM $B_H$ replaced by
TFBM II $B^{I\!I}_{H,\lambda_*}$ and the convergence holds for all $-1/2 <d < \infty $ in contrast to
\cite{sowell2} which is limited to $|d| < 1/2 $.
Under strong tempering $\lambda_N/N \to \lambda_*= \infty $, the limit distribution of
$\widehat \beta_N$ is written in terms of standard Brownian motion but takes a different
form in the cases $d>0, d = 0 $ and $d <0, d \ne \N_- $; moreover, except for the i.i.d. case $d = 0$,
this limit is
different from Sowell's limit in \cite{sowell2} and also
from the unit root distribution
in  Dickey and Fuller~\cite{Dickey1} and Phillips~\cite{Phillips}.

The paper is organized as follows. Section~\ref{sec2} introduces ARTFIMA($p,d,\lambda,q)$ class and provides  basic properties of
these processes. In Section~\ref{sec3} we define TFSM II/TFBM II.
Section~\ref{sec4} contains the main results of the paper (invariance principles).
Section~\ref{sec5} discusses the application to unit root testing.  The proofs of the main results
are relegated to Section \ref{sec6}.

In what follows, $C$ denotes generic constants 
which may be different at different locations.
We write $\limd, \eqd, \limfdd, $  $   \eqfdd $  
for the weak convergence and equality 
of distributions and finite-dimensional
distributions. $\N_\pm := \{\pm 1, \pm 2, \dots \}, \,
\R_+ := (0, \infty), \,
(x)_\pm := \max (\pm x, 0),  x \in \R,  \, \int := \int_\R$. \, $L^p(\R) \, (p \ge 1)$ denotes the Banach space
of measurable functions $f: \R \to \R$ with finite norm $\|f\|_p
= \big(\int |f(x)|^p \d x \big)^{1/p} $.

\section{Tempered fractionally integrated process}\label{sec2}

In this section, we define
ARTFIMA$(p,d,\lambda,q)$ process and discuss its basic properties.
Let $\Phi(z) =1- \sum_{i=1}^ p \phi_i z^i$ and $\Theta(z)=1+ \sum_{i=1}^q \theta_i z^i$ be polynomials
with real coefficients of degree  $p, q\ge 0,$ such that $\Phi (z) $ does not vanish on $\{z \in \C, |z| \le 1 \}$ and
$\Phi(z)$ and $\Theta(z)$ have no common zeros. Let  $d \in \R\setminus \N_+ $.
Consider Taylor's expansion
\begin{equation}
\frac{\Theta(z)}{\Phi(z)}(1-z)^d = \sum_{k=0}^\infty a_{d}(k) z^k , \qquad |z| < 1.
\end{equation}
Note that 
\begin{equation}\label{acoef}
a_{d}(k) = \sum_{s=0}^{k} \omega_{d}(k) \psi(k-s), \qquad k \ge 0, 
\end{equation}
where  $\omega_{d}(k) =  \frac{\Gamma(k-d)}{\Gamma(k+1)\Gamma(d)}$,
and $\psi(j)$ are the coefficients of the power series $\sum_{j=0}^\infty \psi(j) z^j = \Theta(z)/\Phi(z), |z| \le 1 $.
We use the fact (see Kokoszka and Taqqu (\cite{kokoszka}, Lemma~3.1) 
that for any $d \in \R\setminus \N_- $
\begin{equation}\label{omega0}
\omega_{-d}(k) =  \frac{\Gamma(k+d)}{\Gamma(k+1)\Gamma(d)} \ = \ \Gamma(d)^{-1} k^{d-1}\big(1 + O(1/k)\big), \quad k \to \infty. 
\end{equation}

\begin{prop} \label{omegabdd}
Let $\Theta(1) \ne 0$.
For any $d \in \R \setminus \N_- $, the coefficients
$a_{-d}(k), k \ge 0 $ satisfy conditions \eqref{bdk}-\eqref{bdk0}. In particular, for any $d \in \R \setminus \{0,-1,-2,\dots \}$
\begin{equation}\label{omegapdq}
a_{-d}(k) \sim \frac{\Theta(1)}{\Phi(1) \Gamma(d)} k^{d-1}\Big(1 + O\big(\frac{1}{k}\big)\Big) , \qquad k \to \infty.
\end{equation}

\end{prop}

\noi {\it Proof. } Let us  prove \eqref{omegapdq}.
By definition, $a_{-d}(k) = \sum_{j=0}^k
\psi(j) \omega_{-d}(k-j), k \ge 0, $ see \eqref{acoef}. 
It is well-known that $|\psi(j)| \le C \e^{- c j} $ for some constants
$C, c >0$, see (\cite{kokoszka}, proof of Lemma~3.2). Note $\Theta(1)/\Phi(1) = \sum_{j=0}^\infty \psi(j) \ne 0$. We have
\begin{eqnarray*}
&&\big|a_{-d}(k) - (\Theta(1)/\Phi(1)) \omega_{-d}(k) \big| \\
&&= \ \big| \sum_{j=0}^k
\psi(j) \big(\omega_{-d}(k-j)- \omega_{-d}(k)\big) - \omega_{-d}(k) \sum_{j>k} \psi(j) \big|\ \le \ \sum_{i=1}^3
\ell_{k,i},
\end{eqnarray*}
where
\begin{eqnarray*}
\ell_{k,1}&:=&\sum_{0\le j \le k^{1/4}} |\psi(j)| \big|\omega_{-d}(k-j)- \omega_{-d}(k)\big|, \\
\ell_{k,2}
&:=&\sum_{k^{1/4} < j \le k} |\psi(j)| \big( |\omega_{-d}(k-j)| + |\omega_{-d}(k)|\big),
\end{eqnarray*}
and $\ell_{k,3}:= |\omega_{-d}(k)| \sum_{j>k} |\psi(j)| \le C k^{d-1}  \sum_{j>k} \e^{- c j} =  o(k^{d-2}) $ since
$\psi(j)$ decay exponentially.  Similarly,
$\ell_{k,2}
\le $  $  C \e^{- c k^{1/4}} k \max_{1 \le j \le k}   |\omega_{-d}(j)|
=  o(k^{d-2})$.  Using \eqref{omega0} we obtain $|\ell_{k,1}| \le C(\ell'_{k,1} + \ell''_{k,2})$,
where  $\ell''_{k,2} :=  k^{d-2} \sum_{j \ge 0} |\psi(j)| = O(k^{d-2}) $ and
\begin{eqnarray*}
\ell'_{k,1}&:=&\sum_{0\le j \le k^{1/4}} |\psi(j)| \big(k^{d-1} - (k-j)^{d-1}\big)
= k^{d-1} \sum_{0\le j \le k^{1/2}} |\psi(j)| \Big(1 - \big(1 - \frac{j}{k})^{d-1} \Big) \\
&\le&C k^{d-2} \sum_{j\ge 0} j |\psi(j)| =  O(k^{d-2})
\end{eqnarray*}
proving $\big|a_{-d}(k) - (\Theta(1)/\Phi(1)) \omega_{-d}(k) \big| = O(k^{d-2}) $ and hence
\eqref{omegapdq} in view of \eqref{omega0}. Thus, 
$a_{-d}(k)$, $d \ne 0$ satisfy \eqref{bdk}. Condition \eqref{bdk0} is obvious from
$a_{0}(k) = \psi(k)$ and properties of $\psi(k)$ stated above.

It remains to prove \eqref{bdkneg}.  Let $ j <  -d  <   j+1 $ for some $j = 0,1, \ldots. $ Then since
$\Psi(z) := \Theta(z)/\Phi(z)$ is analytic on $\{z \in \C, |z| \le 1 + \delta \}\  (\exists \, \delta >0) $ and
$|a_{-d}(k)| \le C k^{d-1} $, see \eqref{omegapdq}, the function $\sum_{k=0}^\infty a_{-d}(k) z^k
= \Psi (z) (1-  z)^{-d}$ is  $j$ times differentiable on $\{|z| \le 1 \} $ and
$\frac{\partial^i}{\partial z^i }
\Psi (z) (1-  z)^{-d}\big|_{z=1} = \sum_{r=0}^i {i \choose r}
\frac{\partial^r \Psi (z)}{\partial  z^r} \big|_{z=1}
\frac{\partial^{i-r} (1-z)^{-d}}{\partial  z^{i-r}} \big|_{z=1} = 0,  0\le i \le j $. Hence,
$0 = \sum_{k=i}^\infty k (k-1) \cdots (k-i +1) a_{-d}(k)  =
\sum_{k=0}^\infty k (k-1) \cdots (k-i +1) a_{-d}(k)   $ for any $ 0 \le i \le j$, proving
\eqref{bdkneg} and the proposition, too.  \hfill $\Box$

\begin{defn}\label{defn:artfima definition} Let the autoregressive polynomials $\Theta(z), \Phi(z) $ of degree $p,q$  satisfy the above
conditions, and $d \in \R\setminus \N_-,
\lambda \ge 0$. Moreover,
let $\zeta = \{\zeta(t), t \in \Z\} $ be a stationary process with $\E |\zeta(0)| < \infty $.
By  ARTFIMA$(p,d, \lambda,q)$ process with innovation process $\zeta$ we mean a stationary moving-average process
$ X_{p,d,\lambda,q} = \{X_{p,d,\lambda,q}(t), t \in \Z\}$ defined by
\begin{equation}\label{eq:Xjdefinition}
X_{p,d,\lambda,q}(t) = \sum_{k=0}^{\infty}  \e^{-\lambda k}a_{-d}(k)  \zeta(t-k), \qquad t \in \Z
\end{equation}
where the series converges in $L_1 $.

\end{defn}

\begin{rem}
{\rm (i) For $\lambda =0 $ and $|d| < 1/2$  and i.i.d. innovations with zero mean and unit variance,
ARTFIMA$(p,d, 0,q)$ process $X_{p,d,0,q}$ coincides with ARFIMA$(p,d,q)$ process, see e.g.  Brockwell and Davis \cite{BrockwellDavisTSTM}.
Particularly,  
$X_{d,0} = X_{0,d,0,0}$ in \eqref{eq:Xjdefinition}  is a stationary solution of the AR($\infty$) equation
\begin{equation}\label{arfima}
(1- B)^d X_{d,0}(t) = \sum_{j=0}^\infty \omega_{d}(j) X_{d,0}(t-j)
= \zeta(t)
 \end{equation}
and the series  in \eqref{eq:Xjdefinition} and \eqref{arfima} converge in $L_2$, meaning
that $X_{d,0}$ is invertible.

\smallskip

\noi (ii) For $\lambda =0$ and zero mean i.i.d. $\alpha$-stable innovations, $1< \alpha < 2 $,
the definition of $X_{p,d,0,q}$ in Definition \ref{defn:artfima definition}  agrees with the definition of ARFIMA$(p,d,q)$ process
in \cite{kokoszka}, who showed that the series in \eqref{eq:Xjdefinition}
converges a.s. and in $L_1$ for $ d < 1 - \frac{1}{\alpha}$.

 }
\end{rem}

\begin{prop} Let $\Theta(z), \Phi(z) $ satisfy the above
conditions  and $\lambda >0$. Then
the series  in \eqref{eq:Xjdefinition} converges in $L_1 $ for any $d \in \R \setminus \N_-$
hence  ARTFIMA$(p,d, \lambda,q)$ process
$ X_{p,d,\lambda,q} $ in Definition \ref{defn:artfima definition}
is well-defined for arbitrary (stationary) innovation process $\zeta$ with finite mean. Moreover,
if $|\Theta(z)| > 0, \, |z| \le 1 $ then   $ X_{p,d,\lambda,q} $ is invertible:
\begin{equation} \label{Xinv}
\zeta(t) = \sum_{k=0}^\infty \e^{-\lambda k} \widetilde{a_{d}}(k) X_{p,d,\lambda,q}(t-k),
\end{equation}
where
\begin{equation} \label{omegainv}
\sum_{k=0}^\infty \widetilde{a_{d}}(k)\ z^k  = \frac{\Phi(z)}{\Theta(z)}(1 -z)^d, \qquad |z| < 1
\end{equation}
and the series in \eqref{Xinv} converges in $L_1$.

\end{prop}

\noi {\it Proof.} The convergence in $L_1$ of the series in \eqref{eq:Xjdefinition} follows from
\eqref{omegapdq}. To show invertibility of these series, note that by Proposition \ref{omegabdd}
the coefficients in \eqref{omegainv} satisfy the bound $|\widetilde{a_{d}}(k)| \le C k^{-d-1}, k \ge 1 $ for $d \in \R\setminus \N_+ $,
and an exponential bound $ |\widetilde{a_{d}}(k)| \le C \e^{-c k}, k \ge 1, c >0 $ for $d \in  \N_+ $. Hence,
 $\e^{-\lambda k} |\widetilde{a_{d}}(k)| \le C \e^{-\lambda k} k^{-d-1} $, for any $d \in \R$ implying the convergence in $L_1$
 of the series in \eqref{Xinv}. Finally, equality in \eqref{Xinv} follows from identity
 $1 =  \big( \frac{ \Phi(z)}{\Theta(z)}(1 -z)^d \big) \big(\frac{\Theta(z)}{\Phi(z)}(1 -z)^{-d}\big), |z| < 1 $.
 \hfill $\Box$

\smallskip

Proposition \ref{thm:Proprties}  describes some second order properties of ARTFIMA$(p,d,\lambda,q)$ with
standardized innovations.

\begin{prop}\label{thm:Proprties}
{\emph Let $X_{p,d,\lambda,q}$ be ARTFIMA$(p,d,\lambda,q)$ process in \eqref{eq:Xjdefinition}, $d\in \R \setminus \N_-, \lambda>0,$
with standardized i.i.d. innovations $\{\zeta(t), t \in \Z\}, \E \zeta(0) =0, \E \zeta^2(0) = 1$. Then
$\E X_{p,d,\lambda,q}(t) = 0,
\E X^2_{p,d,\lambda,q}(t) = \sum_{k=0}^\infty \e^{-2\lambda k} a^2_{-d}(k) < \infty $ and

\smallskip

\noi (i)  The spectral density of $X_{p,d,\lambda,q}$ is given by
\begin{equation*}
h(x) \ = \  \frac{1}{2\pi}\left|\frac{\Theta(\e^{-\i x})}{\Phi(\e^{-\i x})}\right|^2
(1-2\e^{-\lambda}\cos{x}+ \e^{-2\lambda})^{-d}, \qquad  -\pi\leq x\leq \pi.
\end{equation*}

\noi (ii) The covariance function of $X_{0,d,\lambda,0}$ is given by
\begin{equation}\label{eq:covarianceartfima}
\gamma_{d,\lambda}(k) = \mathbb{E} X_{0,d,\lambda,0}(0) X_{0,d,\lambda,0}(k) = \frac{ \e^{-\lambda k}\Gamma(d+k)}{\Gamma(d)\Gamma(k+1)}\ {_2F_1(d,k+d;k+1;\e^{-2\lambda})},
\end{equation}
where $_2F_1(a, b;c;z)$ is the Gauss hypergeometric function (see e.g. \cite{Gradshteyn}). Moreover,
\begin{equation}\label{eq:shortmemory}
\sum_{k\in \Z}\big|\gamma_{d,\lambda}(k)\big|<\infty, \quad \quad
\sum_{k \in \Z} \gamma_{d,\lambda}(k)= (1-\e^{-\lambda})^{-2d}
\end{equation}
and
\begin{equation}\label{gammalim}
\gamma_{d,\lambda}(k) \sim  A k^{d-1} \e^{-\lambda k}, \quad k \to \infty, \quad \text{where}  \quad
A = (1 - \e^{-2\lambda})^{-d} \Gamma(d)^{-1}.
\end{equation}
}
\end{prop}
{\it Proof.}
(i) From the transfer function $ (\Theta(\e^{-\i x})/\Phi(\e^{-\i x}))  (1- \e^{-\i x - \lambda})^{-d}$
of the filter  in \eqref{eq:Xjdefinition} we have that
$h(x)=\frac{1}{2\pi} |\Theta(\e^{-\i x})/\Phi(\e^{-\i x})|^2  |1-  \e^{-\i x -\lambda }|^{-2d}, $ where
$|1- \e^{-\i x -\lambda}|^2 = 1-2\e^{-\lambda}\cos(x)+ \e^{-2\lambda}$, proving part (i).
(ii) \eqref{eq:covarianceartfima}  follows from $\gamma_{d,\lambda}(k) =\int_{-\pi}^{\pi}\cos{(kx)}h(x)\, \d x $ and
(\cite{Gradshteyn}, Eq.\ 9.112). The first relation in \eqref{eq:shortmemory} follows from
$\sum_{j \in \Z} |e^{-\lambda j} \omega_{d}(j)| < \infty $, see \eqref{omega0}, and the second
one from $\sum_{k\in \Z}\gamma_{d,\lambda}(k)= 2\pi h(0)=(1- \e^{-\lambda})^{-2d}$.
Finally, \eqref{gammalim} is proved in (\cite{giraitis}, (4.15)).
\hfill $\Box$

\section{Tempered fractional Brownian and stable motions of second kind}\label{sec3}

This section contains the definition of  TFBM II/TFSM II and some of its properties
from Sabzikar and Surgailis~\cite{SurgailisFarzadTFSMII}. The reader is referred to the aforementioned paper
for further properties
of these processes including relation to tempered fractional calculus, relation between
TFBM II/TFSM II and TFBM/TFSM, dependence properties of the increment process
(tempered fractional Brownian/stable noise), local and global asymptotic self-similarity.

For $1< \alpha \le 2$,
let  $M_\alpha = \{M_\alpha(t), t \in \R\}$ be an $\alpha$-stable 
L\'evy process with stationary independent increments and characteristic
function
\begin{equation}\label{Mcf}
\E \e^{\i \theta M_\alpha (t)}  = \e^{ - \sigma^\alpha |\theta|^\alpha |t|
\big(1 - \i \beta \tan (\pi \alpha/2) {\rm sign}(\theta)\big)},
\qquad \theta \in \R,
\end{equation}
where $\sigma >0 $ and $ \beta \in [-1,1]$
are  the scale and skewness parameters, respectively. For $\alpha =2$,
$M_2(t) = \sqrt{2} \sigma B(t)$, where $B $ is a standard Brownian motion with variance $\E B^2(t) = t$. Stochastic
integral $I_\alpha(f) \equiv \int f(x) M_\alpha (\d x)$ is defined for any $ f \in L^\alpha (\R)$
as $\alpha$-stable
random variable with
characteristic function
\begin{equation} \label{Ialpha}
\E \e^{ \i \theta I_\alpha (f)}  = \exp \{  - \sigma^\alpha |\theta|^\alpha
\int |f(x)|^\alpha
\big(1 - \i \beta \tan (\pi \alpha/2) {\rm sign}(\theta f(x))\big)  \d x \},
\quad \theta \in \R.
\end{equation}
see e.g. \cite[Chapter 3]{SamorodnitskyTaqqu}.  \\
For $t \in \R, H> 0, 1 < \alpha \le 2,
\lambda \ge 0$ consider
the function $y \mapsto h_{H,\alpha,\lambda}(t;y): \R \to \R $ given in \eqref{hdef0}.
Note $h_{H, \alpha,\lambda}(t;\cdot) \in L^\alpha(\R)  $ for any $t \in \R, \lambda >0, 1< \alpha \le 2,  H > 0 $  and also for
$\lambda = 0, 1< \alpha \le 2,
H \in (0,1)$.
We will use the following integral representation of \eqref{hdef0} (see \cite{SurgailisFarzadTFSMII}). \\
For $H > \frac{1}{\alpha}$:
\begin{eqnarray}\label{hdef1}
&h_{H,\alpha,\lambda}(t;y) = (H - \frac{1}{\alpha})
\int_0^t (s-y)_+^{H - \frac{1}{\alpha} -1} \e^{-\lambda (s-y)_+} \d s.
\end{eqnarray}
For $0< H < \frac{1}{\alpha}$:
\begin{eqnarray} \label{hdef2}
&h_{H,\alpha,\lambda}(t;y)=(H - \frac{1}{\alpha})
\begin{cases}\int_0^t (s-y)_+^{H - \frac{1}{\alpha} -1} \e^{-\lambda (s-y)_+} \d s, &y < 0, \\
-\int_t^\infty (s-y)_+^{H - \frac{1}{\alpha} -1} \e^{-\lambda (s-y)_+} \d s + \lambda^{\frac{1}{\alpha} -H}
\Gamma(H - \frac{1}{\alpha})
, &y \ge 0.
\end{cases}
\end{eqnarray}

\begin{defn}\label{defTHP}
{\rm
{Let $M_\alpha $ be $\alpha$-stable L\'evy process in \eqref{Mcf}, $1< \alpha \le 2 $ and
$H >0, \, \lambda > 0$.
Define, for $t\in \R$, the stochastic integral
\begin{eqnarray}\label{TFSM2}
Z^{I\!I}_{H,\alpha,\lambda}(t)&:=&\int h_{H,\alpha,\lambda}(t;y)
M_\alpha (\d y).
\end{eqnarray}
The process $Z^{I\!I}_{H,\alpha,\lambda}  = \{ Z^{I\!I}_{H,\alpha,\lambda}(t), t\in \R\}$ will be called
{\it tempered fractional stable motion of  second kind} (TFSM II).
A particular case of \eqref{TFSM2} corresponding to $\alpha = 2$
\begin{eqnarray}\label{TFBM2}
B^{I\!I}_{H,\lambda}(t)
&:=&\frac{1}{\Gamma(H+\frac{1}{2})}\int  h_{H,2,\lambda}(t;y) B(\d y)
\end{eqnarray}
will be called  {\it tempered  fractional Brownian motion of second kind} (TFBM II). } }
\end{defn}

The next proposition states some basic properties of TFSM II \, $Z^{I\!I}_{H,\alpha,\lambda}$.

\begin{prop} {\rm (\cite{SurgailisFarzadTFSMII})} \label{lem:g squar integrable}

\smallskip

\noi (i) $Z^{I\!I}_{H,\alpha,\lambda}$
in \eqref{TFSM2} is well-defined for any $t \in \R $ and
$1< \alpha \le 2, H >0, \lambda > 0$,
as a stochastic integral in \eqref{Ialpha}.

\smallskip

\noi (ii) $Z^{I\!I}_{H,\alpha,\lambda}$
in \eqref{TFSM2} has stationary increments and
$\alpha$-stable finite-dimensional distributions. Moreover, it satisfies the following
scaling property:
$\big\{Z^{I\!I}_{H, \alpha, \lambda}(bt), t \in \R\big\} \eqfdd
\big\{b^{H}Z^{I\!I}_{H,\alpha, b\lambda}(t),  t \in \R\big\}, $  $ \forall \ b >0. $

\smallskip

\noi (iii) $Z^{I\!I}_{H,\alpha,\lambda}$
in \eqref{TFSM2} has
a.s. continuous paths if either $\alpha = 2,  H>0$, or $ 1< \alpha < 2,
H > 1/\alpha$ hold.



\end{prop}

\section{Invariance principles }\label{sec4}

In this section, we discuss invariance principles and  the convergence of (normalized) partial sums process
\begin{equation}\label{eq:S definition}
S^{d,\lambda}_{N}(t):=\sum_{k=1}^{[Nt]}X_{d,\lambda}(k),  \qquad t\in [0,1]
\end{equation}
of tempered linear process  $X_{d,\lambda}$ in \eqref{TL}  with i.i.d. innovations belonging to the domain of attraction
of $\alpha$-stable law, $1< \alpha \le 2 $ (see below). We shall assume that the tempering parameter
$\lambda \equiv \lambda_N \to  0$ as $N \to \infty$
and following limit exists:
\begin{equation}\label{lambdaC}
\lim_{N \to \infty}  N \lambda_N = \lambda_*  \in [0, \infty].
\end{equation}
Recall that
$X_{d,\lambda_N}$ is said
{\it strongly tempered} if $\lambda_* =  \infty$,
{\it weakly tempered} if $\lambda_* =  0$, and
{\it moderately tempered} if $0 < \lambda_* < \infty $.  We show that the limits of
the partial sums process in \eqref{eq:S definition} exist under condition \eqref{lambdaC} and depend on $\lambda_*,
d, \alpha$; moreover, in all cases these limits belong to the class of TFSM II
processes defined in Section~3.

\begin{defn}\label{Zdef} Write $\zeta \in D(\alpha)$, $1< \alpha \le 2$ if

\begin{description}

\item[{\rm (i)}]  $\alpha = 2$ and $\E \zeta = 0, \  \sigma^2 := \E \zeta^2 <
\infty$, or
\item[{\rm (ii)}] $1 < \alpha < 2 $ and there exist some constants
$c_{1}, c_{2} \ge 0, c_{1} + c_{2} > 0$
such that $\lim_{x \to \infty} x^\alpha {\rm P}(\zeta > x) =c_{1} $ and $\lim_{x \to -\infty} |x|^{\alpha}
{\rm P}(\zeta \le x)= c_{2};$
moreover, $\E \zeta = 0$.
\end{description}
\end{defn}

Condition $ \zeta \in D(\alpha$) implies that  r.v. $\zeta$ belongs to the domain of normal attraction of an
$\alpha$-stable law. In other words, if $\zeta(i), i \in \Z$ are i.i.d. copies of $\zeta $ then
\begin{equation}\label{clt}
N^{-1/\alpha}\sum_{i=1}^{[Nt]} \zeta(i)
\ \limfdd  \ M_\alpha(t), \qquad N \to \infty,
\end{equation}
where $M_\alpha$ is an $\alpha$-stable L\'evy process  in \eqref{Mcf} with  $\sigma, \beta  $
determined by $c_1, c_2 $, see (\cite{fell1966}, pp.~574-581).
We shall  use the following criterion
for convergence of weighted sums in i.i.d. r.v.s. See (\cite{koul}, Prop.~14.3.2), \cite{astrauskas}, \cite{kasahara}.

\begin{prop} \label{propQ} Let $1< \alpha \le 2 $ and
$Q(g_N) = \sum_{t \in \Z} g_N(t) \zeta(t) $ be a linear form in i.i.d. r.v.s
$\zeta(t) \in D(\alpha)$ with real coefficients $g_N(t), t \in \Z$. Assume that
there exists $p \in [1,\alpha) $ if $\alpha<2$, $p=2 $ if $\alpha =2 $ and
a function $g \in L^p (\R) $ such that the functions
\begin{equation}\label{gNx}
\widetilde g_N(x) := N^{1/\alpha} g_N([xN]), \qquad x \in \R
\end{equation}
satisfy
\begin{equation}\label{pcond}
\| \widetilde g_N - g \|_p \to 0 \qquad (N \to \infty).
\end{equation}
Then $ Q(g_N) \limd \int g(x) M_\alpha (\d x) $, where $M_\alpha $ is as
in \eqref{clt} and
\eqref{Mcf}.

\end{prop}

Write   $\DD $ for weak convergence of random processes in the Skorohod space
$D[0,1]$ equipped with $J_1$-topology, see \cite{Bill}.

\smallskip

In Theorem \ref{thm:theorem one} below,  $X_{d,\lambda_N}$ is a tempered linear process of  \eqref{TL} with i.i.d.
innovations $\zeta(t) \in D(\alpha), 1 < \alpha \le 2$,
coefficients $b_d(k), k \ge 0, \, d \in \R \setminus \N_-$ satisfying \eqref{bdk}-\eqref{bdk0}, and tempering
parameter $  0< \lambda_N \to 0 \, (N \to \infty)$
satisfying \eqref{lambdaC}. W.l.g., we shall assume
that the asymptotic constant $c_d$ in \eqref{bdk}-\eqref{bdk0} equals 1: $c_d = 1 \ \forall d \in \R \setminus \N_- $.

\begin{thm}\label{thm:theorem one}

\noi (i) {\rm (Strongly tempered process.)} Let $\lambda_* = \infty$ and  $d \in \R \setminus \N_-$. Then
\begin{equation}\label{eq:theorem one first part}
N^{-\frac{1}{\alpha}} \lambda_N^d S^{d,\lambda_N}_{N}(t)\ \limfdd \ M_\alpha(t),
\end{equation}
where $M_\alpha$ is $\alpha$-stable L\'evy process  in \eqref{clt}. Moreover, if $\alpha =2$ and
$\mathbb{E} |\zeta(0)|^{p} <\infty$ for some $p>2$ then
\begin{equation}\label{eq:theorem one second part}
N^{-\frac{1}{2}} \lambda_N^d S^{d,\lambda_N}_{N}(t) \ \DD \ \sigma B(t),
\end{equation}
where $B$ is a standard Brownian motion and $\sigma >0$ some constant.

\smallskip

\noi (ii)  {\rm (Weakly tempered process.)} Let $ \lambda_* = 0 $ and $H = d + \frac{1}{\alpha} \in (0,1)  $. Then
\begin{equation}\label{eq:theorem one 3 part}
N^{-H} S^{d,\lambda_N}_{N}(t) \ \limfdd  \ \text{\small $\Gamma (d+1)^{-1}$} Z_{H,\alpha,0}(t),
\end{equation}
where $Z_{H,\alpha,0}$ is a linear fractional $\alpha$-stable motion, see
Definition~\ref{defTHP}. 
Particularly,
for $\alpha =2 $, $Z_{H,2,0}$ is a multiple of FBM. \\
Moreover, if either
$1 < \alpha \le  2,
1/\alpha < H < 1 $, or  $\alpha =2, 0 < H < 1/2 $ and $\E |\zeta(0)|^{p} < \infty \, (\exists p > 1/H)$  hold, then
$\limfdd$ in \eqref{eq:theorem one 3 part}
 can be replaced by $\DD $.

\smallskip

\noi (iii)  {\rm (Moderately tempered process.)}  Let $ \lambda_* \in (0,\infty) $
and $H = d + \frac{1}{\alpha} >0  $.
Then
\begin{equation}\label{eq:theorem one 4 part}
N^{-H}S^{d,\lambda_N}_{N}(t) \ \limfdd \ \text{\small $\Gamma (d+1)^{-1}$}Z^{I\!I}_{H,\alpha,\lambda_*}(t),
\end{equation}
where $Z^{I\!I}_{H,\alpha,\lambda_*}$ is a TFSM II as defined in
Definition~\ref{defTHP}. \\
Moreover, if either
$1 < \alpha \le  2, 1/\alpha < H $, or  $\alpha =2, 0 < H < 1/2 $ and $\E |\zeta(0)|^{p} < \infty \, (\exists p > 1/H)$  hold, then
$\limfdd$ in \eqref{eq:theorem one 4 part}
can be replaced by $\DD $.
\end{thm}

\begin{rem} {\rm Note that for $\lambda_N = \lambda_*/N $ the normalization
in \eqref{eq:theorem one first part} becomes
$N^{-(\frac{1}{\alpha} +d)}  \lambda_*^d $ where the exponent $\frac{1}{\alpha} + d = H $ is the same as
in \eqref{eq:theorem one 3 part} and
\eqref{eq:theorem one 4 part}.

}
\end{rem}

\begin{rem} {\rm The functional convergence in \eqref{eq:theorem one first part}, case $1< \alpha < 2 $
(the case of discontinuous limit process) is open and
apparently does not hold in the usual $J_1$-topology, see  \cite{balan}.
In the case of \eqref{eq:theorem one 3 part} and \eqref{eq:theorem one 4 part} and $1 < \alpha < 2,  0< H < \frac{1}{\alpha} $,
functional convergence
cannot hold in principle since the limit processes do not belong to $D[0,1]$.

}
\end{rem}

\section{Tempered fractional unit root distribution}\label{sec5}

A fundamental problem of time series is testing for the unit root $\beta = 1$
in the AR(1) model
\begin{equation}\label{AR1}
Y(t) = \beta Y (t-1)+ X(t), \quad\quad t=1,2,\ldots,N, \ \  Y(0) = 0
\end{equation}
with stationary error process $X = \{ X(t), t \in \Z \} $. The classical approach to the unit root testing
is based on the limit distribution of
the OLS estimator $\widehat \beta_N $
\begin{equation}\label{betadef}
\widehat \beta_N = \frac{ \sum_{t=1}^{N} Y (t) Y(t-1)}{\sum_{t=1}^{N} Y^{2}(t-1)}.
\end{equation}
The limit theory for $\widehat \beta_N $ in the case of weakly dependent errors $X$ was developed
in Phillips~\cite{Phillips}. We note that \cite{Phillips} 
makes an extensive use of invariance principle for the error process.
Sowell~\cite{sowell2} obtained the limit distribution of $\widehat \beta_N$
in the case of strongly dependent ARFIMA$(0,d,0)$
error process
with finite variance and  standardized i.i.d. innovations.
\cite{sowell2} proved that the distribution of the normalized statistic
$N^{1 \wedge (1+2d)} (\widehat \beta_N - 1) $ tends to
that of the ratio
\begin{equation} \label{sowell}
\frac{1}{2\int_0^1 B^2_H(s) \d s}
\begin{cases}
B^2_H(1), &0< d < 1/2, \\
B^2(1) - 1, &d=0, \\
- H \Gamma(H + \frac{1}{2})/\Gamma(\frac{3}{2}-H),  &-1/2 < d < 0,
\end{cases}
\end{equation}
where
$H = d + \frac{1}{2}$ and $B_H$ is a  FBM with parameter $H \in (0,1), \, B = B_{1/2}$ being a standard
Brownian motion.

In this section we extend Sowell's \cite{sowell2} result  to
ARTFIMA$(0,d,\lambda_N,0)$ error process with small tempering
parameter $\lambda_N  \sim \lambda_*/N \to 0 $ as in \eqref{lambdaC}.
Although Theorem~\ref{thm:ADDIDTIVE ARTFIMA TO REGRESSION MODEL1} can be generalized
to more general tempered processes with finite variance as in Theorem~\ref{thm:theorem one},
our choice of  ARTFIMA$(0,d,\lambda_N,0)$ as the error process is motivated by
better comparison to \cite{sowell2}.  As noted in Section~1,
the degree of tempering has a strong effect on the limit distribution of
$\widehat \beta_N$ and leads to a new two-parameter family of
tempered fractional unit root distributions. Following \cite{sowell2},
we decompose
\begin{equation}\label{betaN}
\widehat \beta_N -1 = \widehat A_N - \widehat B_N,
\end{equation}
where
\begin{eqnarray}\label{AB}
\widehat A_N := \frac{Y^2(N) }{ 2\sum_{t=1}^{N} Y^{2}(t-1)},
\qquad \widehat B_N := \frac{\sum_{t=1}^N X^2(t) }{ 2\sum_{t=1}^{N} Y^{2}(t-1)}.
\end{eqnarray}
Under the unit root hypothesis $\beta = 1$ we have $Y(t) = \sum_{i=1}^t X(i) = S_N(t/N) $, where
$S_N(x) := \sum_{t=1}^{[Nx]} X(t), x \in [0,1]$ is the partial sums process. Particularly,
the statistics in \eqref{AB} can be rewritten as
\begin{eqnarray}\label{AB1}
\widehat A_N = \frac{S^2_N(1)}{2N \int_0^1 S^2_N(s) \d s}, \qquad
\widehat B_N =  \frac{\sum_{t=1}^N X^2(t)}{ 2N \int_0^1 S^2_N(s) \d s}.
\end{eqnarray}
For ARTFIMA error process
$X_{d,\lambda_N}$, the behavior of $S^2_N(1)$  and $\int_0^1 S^2_N(s) \d s $ can be derived from
Theorem~\ref{thm:theorem one}. The behavior of  $\sum_{t=1}^N X^2(t)$ is established
in the following proposition.

\begin{prop} \label{propX2} Let $X_{d,\lambda_N} $ be an ARTFIMA$(0,d, \lambda_N,0)$ process in \eqref{eq:Xjdefinition} with i.i.d.
innovations $\{\zeta (t)\}, \E \zeta(0) = 0, \E \zeta^2(0) = 1$,
fractional parameter $d \in \R \setminus \N_-$ and tempering
parameter $\lambda_N \to 0$.  Moreover, let $\E |\zeta(0)|^p < \infty \, (\exists \, p > 2)$. 
Then
\begin{eqnarray}\label{sum1}
\frac{1}{N}\sum_{t=1}^N X^2_{d,\lambda_N}(t)
&\limp&\frac{\Gamma(1-2d)}{\Gamma^2(1-d)},  \qquad  d < 1/2, \\
\label{sum2}
\frac{\lambda_N^{2d-1}}{N}
\sum_{t=1}^N X^2_{d,\lambda_N}(t)&\limp&
\frac{\Gamma(d-1/2)}{2\sqrt{\pi}\,\Gamma(d)},  \qquad d > 1/2, \\
\label{sum3}
\frac{1}{N |\log \lambda_N|}
\sum_{t=1}^N X^2_{d,\lambda_N}(t)&\limp&\frac{1}{\pi}, \hskip2cm d=1/2.
\end{eqnarray}

\end{prop}

\noi The main result of this section is the following theorem.

\begin{thm}\label{thm:ADDIDTIVE ARTFIMA TO REGRESSION MODEL1}
Consider the AR(1) model in \eqref{AR1} with $\beta = 1$ and ARTFIMA$(0,d,\lambda_N,0)$ error process
$X = \{X_{d,\lambda_N}(t)\} $ in \eqref{eq:Xjdefinition} with i.i.d.
innovations $\{\zeta(t), t \in \Z\}, \E \zeta(0) = 0, \E \zeta^2(0) = 1,  \E |\zeta(0)|^p < \infty \ (\exists \, p > 2 \vee 1/(d+1/2))$,  fractional parameter $d \in \R \setminus \N_-$ and tempering
parameter $\lambda_N >0$
satisfying \eqref{lambdaC}.

\medskip

\noi (i) {\rm (Strongly tempered errors.)} Let $\lambda_* = \infty, \, d \in \R  \setminus \N_- $.  Then
\begin{eqnarray*}\label{eq:unit root part I}
\min(1, \lambda^{-2d}_N)   N(\widehat\beta_N-1)
&\limd&\frac{1}{2\int_{0}^{1}B^2(s)\, \d s}
\begin{cases}
B^2(1),  &d > 0, \\
B^2(1)-1, &d=0, \\
-\Gamma(1-2d)/\Gamma(1-d)^2, &d<0,
\end{cases}
\end{eqnarray*}
where $B$ is a standard Brownian motion.

\smallskip

\noi (ii) {\rm (Weakly tempered errors.)} Let $ \lambda_* = 0 $ and $H = d+ \frac{1}{2}\in (0,1)$.  Then
\begin{eqnarray}\label{eq:unit root part II}
N^{1 \wedge (1+2d)}(\hat{\beta}_N -1)&\limd&\frac{1}{2\int_0^1 B^2_H(s) \d s }
\begin{cases}
B^2_H(1), &\frac{1}{2}< H <1, \\
B^2(1)-1,&  H=\frac{1}{2} \\
- H \Gamma(H + \frac{1}{2})/\Gamma(\frac{3}{2}-H), &0<H<\frac{1}{2},
\end{cases}
\end{eqnarray}
where $B_{H}$ is a FBM with variance $\E B^2_H(t) = t^{2H}, B = B_{1/2}$.

\smallskip

\noi (iii) {\rm (Moderately tempered errors.)} Let $ 0<\lambda_* <\infty$ and
$H = d+ \frac{1}{2} >0$. Then
\begin{eqnarray}\label{eq:unit root part III}
N^{1 \wedge (1+2d)}(\hat{\beta}_N-1)&\limd&\frac{1}{2\int_0^1 \big(B^{I\!I}_{H,\lambda_*}(s)\big)^2 \d s }
\begin{cases}
(B^{I\!I}_{H,\lambda_*}(1))^2, &H> \frac{1}{2}, \\
(B^{I\!I}_{\frac{1}{2},\lambda_*}(1))^{2}-1,  &  H=\frac{1}{2} \\
-\Gamma(2(1-H))/\Gamma(\frac{3}{2}-H)^2, & 0< H < \frac{1}{2},
\end{cases}
\end{eqnarray}
where $B^{I\!I}_{H,\lambda}$ is a TFBM II given by \eqref{TFBM2}.

\end{thm}

\begin{rem}
{\rm The limit  \eqref{eq:unit root part II} in the weakly tempered case
coincides with Sowell's
limit \eqref{sowell}.  Since $ B^{I\!I}_{H,0} \eqfdd C_1 B_H $, where
$C_1^2 =  \Gamma (1-H)\big/ 2^{2H} H \Gamma (H+ 1/2) \sqrt{\pi}$,
see \cite{SurgailisFarzadTFSMII},
the  r.v. on the r.h.s. of \eqref{eq:unit root part III} for $\lambda_* =0, 0< H< 1 $
also coincides with   \eqref{sowell}, however the convergence  \eqref{eq:unit root part III}
holds for any $H >0$ in contrast to $H\in (0,1)$  in  \eqref{eq:unit root part II}.

}
\end{rem}

\noi {\it Proof of Theorem \ref{thm:ADDIDTIVE ARTFIMA TO REGRESSION MODEL1}.} From Theorem~\ref{thm:theorem one} and
Proposition~\ref{propX2} we obtain the joint convergence of
\begin{eqnarray}\label{joint}
&\Big(a_N (S^{d,\lambda_N}_N(1))^2, \,
a_N \int_0^1 (S^{d,\lambda_N}_N(s))^2 \d s, \, b_N \sum_{t=1}^N X^2_{d,\lambda_N}(t)\Big)
\end{eqnarray}
where $a_N \to 0, b_N \to 0$ are normalizations defined in these theorems and
depending on $d $ and $\lambda_*$, in each case
(i)-(iii) of Theorem~\ref{thm:ADDIDTIVE ARTFIMA TO REGRESSION MODEL1}. Then the statement of
Theorem~\ref{thm:ADDIDTIVE ARTFIMA TO REGRESSION MODEL1} follows from \eqref{joint}, the continuous mapping
theorem   and the representation of $\widehat \beta_N -1 $ in \eqref{betaN}-\eqref{AB1}
through the corresponding quantities in \eqref{joint}.
\hfill $\Box$

\section{Proofs of Theorem~\ref{thm:theorem one} and Proposition~\ref{propX2}
}\label{sec6}

\noi {\it Proof of Theorem~\ref{thm:theorem one}.} (i) We restrict the proof of finite-dimensional convergence
in \eqref{eq:theorem one first part} to one-dimensional convergence at $t >0$
since the general case follows similarly.  We use Proposition~\ref{propQ}. Accordingly, write
$ N^{-\frac{1}{\alpha}} \lambda_N^d S^{d,\lambda_N}_{N}(t) = Q(g_N(t,\cdot)) = \sum_{i \in \Z} g_N(t;i) \zeta(i),$ where
$g_N(t;i) := $  $ N^{-1/\alpha} \lambda_N^d $  $ \sum_{k=1\vee i}^{[Nt]} \e^{-\lambda_N (k-i)}\ b_d(k-i) $. It suffices to prove
\eqref{pcond} for
suitable $p$ and $g(t;x) := \1_{[0,t]}(x)$.
We have
\begin{equation}
\widetilde g_N(t;x) = \lambda_N^d \sum_{k=1 \vee [Nx]}^{[Nt]} \e^{-\lambda_N (k - [Nx])}\ b_d(k- [Nx]) .
\end{equation}
Let us prove the point-wise convergence:
\begin{equation}\label{gNlim}
\widetilde g_N(t;x) \to g(t;x) =  \1_{[0,t]}(x), \quad \forall \  x \ne 0, t.
\end{equation}
Let us prove that conditions \eqref{bdk}-\eqref{bdk0} imply that
\begin{eqnarray}\label{GNlim}
&&G_N:=\lambda_N^{d} \sum_{k=0}^\infty \e^{-\lambda_N k}\,  b_d(k)
\to  1  \qquad (N \to \infty).   
\end{eqnarray}
First, let $d >0$. 
Then since $\sum_{k=0}^n b_d(k) 
\sim (1/d \Gamma (d)) n^d, n \to \infty $ according to \eqref{bdk},
then by applying the Tauberian theorem for power series (Feller \cite{fell1966}, Ch.~13, \S~5, Thm.~5) we have
$\sum_{k=0}^{\infty}b_d(k) \,\e^{-\lambda_N k}\sim ( 1- \e^{-\lambda_N} )^{-d},  \, N\to\infty$,
proving \eqref{GNlim} for $d>0$.
Next, let $d=0$. Then in view of \eqref{bdk0} the dominated convergence theorem applies yielding
$\sum_{k=0}^\infty \e^{-\lambda_N k}\,  b_0(k) \to  \sum_{k=0}^\infty b_0(k) = 1$ and
\eqref{GNlim} follows again.

Next, let $-1 < d < 0$. Then $\tilde b_d(k) := \sum_{i=k}^\infty b_d(i) \sim (-1/d \Gamma (d)) k^{\tilde d -1}, \, k \to \infty,
\tilde d := d+1 \in  (0,1), \,  \tilde b_d(0) = 0$ and
$\sum_{k=0}^\infty \e^{-\lambda k}  b_d(k) = -\e^{\lambda} (1- \e^{-\lambda}) \sum_{k=1}^\infty
\e^{-\lambda k} \tilde b_d(k) $ using summation by parts. Then the aforementioned Tauberian theorem
implies $\sum_{k=0}^\infty \e^{-\lambda_N k}  b_d(k) \sim   (1- \e^{-\lambda_N})^{1 - \tilde d}
  =   (1- \e^{-\lambda_N})^{-d}  $ proving \eqref{GNlim} for
$-1 < d < 0$. In the general case
$-j < d < -j +1, j =1,2, \dots $ relation
\eqref{GNlim} follows similarly using
summation by parts $j$ times.


Let  $ 0< x < t$ first.  Then $\widetilde g_N(t;x) = G_N - \widetilde g^*_N(t;x), $ where
\begin{eqnarray*}
\widetilde g^*_N(t;x) := \lambda_N^d \sum_{k> [Nt]- [Nx]} b_d(k) \e^{-\lambda_N k}.
\end{eqnarray*}
Using \eqref{bdk} for $d \ne 0$  we obtain
\begin{eqnarray}
|\widetilde g^*_N(t;x)|&\le&C (N\lambda_N)^d N^{-1} \sum_{k> [Nt]- [Nx]}  \e^{- (N\lambda_N) (k/N)} (k/N)^{d -1} \nn \\
&\le&C(N\lambda_N)^d \int_{t-x}^\infty \e^{- N \lambda_N y } y ^{d-1} \d y  \nn \\
&=&C\int_{(t-x)(N \lambda_N) }^\infty \e^{- z } z^{d-1} \d z  \to  0 \label{gN1}
\end{eqnarray}
since $ N \lambda_N \to \infty $. A similar result for $d =0$ follows directly from \eqref{bdk0}.
In view of \eqref{GNlim}, this proves \eqref{gNlim} for $0< x <t$.
Next, let $x < 0$.  Then similarly as above
\begin{eqnarray}
\widetilde g_N(t;x)&\le&C (N\lambda_N)^d N^{-1} \sum_{k> |[Nx]|}  \e^{- (N\lambda_N) (k/N)} (k/N)^{d -1} \nn \\
&\le&C(N\lambda_N)^d \int_{|x|}^\infty \e^{- N \lambda_N y } y ^{d-1} \d y \nn \\
&=&C\int_{|x|(N \lambda_N) }^\infty \e^{- z } z^{d-1} \d z  \to  0, \label{gN2}
\end{eqnarray}
proving \eqref{gNlim}.  Note also that
$|\widetilde g_N(t;x)| \le C \lambda_N^d \sum_{k=0}^\infty |b_d(k)|  \le  C \lambda_N^d \le C $
for $d <0$  and $ |\widetilde g_N(t;x)| \le C (N\lambda_N)^d  \int_{0}^\infty \e^{- N \lambda_N y } $ $ y ^{d-1} \d y
\le C \int_{0}^\infty \e^{- z }  z^{d-1} \d z \le C $ for $d >0$, implying that
$|\widetilde g_N(t;x)|$ is bounded uniformly in $x \in \R, N \ge 1; $
moreover, according to  \eqref{gN2}
$|\widetilde g_N(t;x)| \le C \e^{-c' |x|},  x < -2$  decays exponentially with $x \to  - \infty $
with some $c'>0$  uniformly in $N\ge 1$. This proves
\eqref{pcond} and hence \eqref{eq:theorem one first part}.

Consider the functional convergence in \eqref{eq:theorem one second part}. This follows
from the tightness criterion
\begin{equation}\label{tight1}
N^{-p/2} \lambda^{pd}_N \E |S^{d,\lambda_N}_{N}(t) - S^{d,\lambda_N}_{N}(s)|^p
\le C |L_N(t) - L_N(s)|^{p/2}, \qquad \forall \  0 \le s < t \le 1,
\end{equation}
where $L_N(t) := [Nt]/N$, see  \cite{Bill}, also (\cite{koul}, Lemma~4.4.1). By Rosenthal's inequality (see
e.g. \cite{koul}, Proposition~4.4.3),
$ \E |S^{d,\lambda_N}_{N}(t) - S^{d,\lambda_N}_{N}(s)|^p \le C \E^{p/2} |S^{d,\lambda_N}_{N}(t) - S^{d,\lambda_N}_{N}(s)|^2 $
where $\E |S^{d,\lambda_N}_{N}(t) - S^{d,\lambda_N}_{N}(s)|^2 =   N \lambda_N^{-2d}  \int |\widetilde g_N(t;x) -
\widetilde g_N(s;x)|^2 \d x \le C N \lambda_N^{-2d} |L_N(t) - L_N(s)|   $ follows similarly as above,
proving  \eqref{tight1} and part (i), too.

\smallskip

\noi (ii) Relation \eqref{eq:theorem one 3 part} is well-known with $S^{d,\lambda_N}_{N}(t)$ replaced by
$S^{d,0}_{N}(t)$, see e.g. \cite{astrauskas}, also (\cite{koul}, Cor.~4.4.1),
so that it suffices
to prove
\begin{equation}\label{Sapp}
R_N(t) := S^{d,\lambda_N}_{N}(t) - S^{d,0}_{N}(t) = o_p (N^H).
\end{equation}
With Proposition~\ref{propQ} in mind, \eqref{Sapp} follows from
$\|\widetilde g^0_N(t;\cdot) \|_p \to 0$, where
\begin{eqnarray*} 
|\widetilde g^0_N(t;x)|
&:=&N^{-d} \Big|\sum_{k=1 \vee [Nx]}^{[Nt]} b_d (k- [Nx])\ \big(1- \e^{-\lambda_N (k - [Nx])}\big)\Big|  \\
&\le&C(N \lambda_N) \frac{1}{N^{d+1}}  \sum_{k=1 \vee [Nx]}^{[Nt]} (k - [Nx])^d \le  C (N \lambda_N) \to  0 \nn
\end{eqnarray*}
uniformly in $x \in \R$,
where we used \eqref{bdk}, inequality $1 - \e^{-x} \le x \, (x \ge 0) $ and the fact that
$H \in (0,1), 1 < \alpha \le 2 $ imply $-1 <   d < 1 - \frac{1}{\alpha}  $. For $ x < -1 $ a similar
argument leads to
\begin{eqnarray} \label{g02}
|\widetilde g^0_N(t;x)|
&\le&C N^{-d} \sum_{k=1}^{[Nt]} (k- [Nx])^{d-1}  \ \le C N^{-d} \big(  ([Nt]- [Nx])^{d} - (- [Nx])^d \big) \\
&\le&C ((t-x)^d - (-x)^d) \le C (-x)^{d-1} \nn
\end{eqnarray}
implying the dominating bound $|\widetilde g^0_N(t;x)| \le C/(1 + |x|)^{1-d} =: \bar g(x) $ where
$\| \bar g \|_p  < \infty $ for $1 \le p < \alpha $ sufficiently close to $\alpha $ due to condition
$ d < 1 - \frac{1}{\alpha} $. This proves $\|\widetilde g^0_N(t;\cdot) \|_p \to 0$, hence
\eqref{Sapp} and \eqref{eq:theorem one 3 part}, too.

To prove the tightness part of (ii), we use a similar criterion as in \eqref{tight1}, viz.,
\begin{equation}\label{tight2}
N^{-p H} \E |S^{d,\lambda_N}_{N}(t) - S^{d,\lambda_N}_{N}(s)|^p
\le C |L_N(t) - L_N(s)|^{q}, \qquad \forall \  0 \le s < t \le 1,
\end{equation}
with $L_N(t) = [Nt]/N$ and suitable $p, q > 1 $.
Let first $  \frac{1}{\alpha} <   H < 1,  $ or  $0< d < 1-\frac{1}{\alpha}$.
Let  $\widetilde g_N(t;x) = N^{-d} \sum_{k=1 \vee [Nx]}^{[Nt]} b_d (k- [Nx]) $  $\e^{-\lambda_N (k - [Nx])} $.
Then for $ 0< s < t$
\begin{eqnarray}
|\widetilde g_N(t;x) - \widetilde g_N(s;x)|
&\le&N^{-d} \sum_{k= [Ns] +1}^{[Nt]} |b_d (k- [Nx])|
\le CN^{-d} \sum_{k= [Ns] +1}^{[Nt]} (k- [Nx])_+^{d-1} \nn \\
&\le&CN^{-d} \int_{[Ns]}^{[Nt]} (y-[Nx])_+^{d-1} \d y \nn \\
&\le&C((L_N(t)-L_N(x))_+^d  - (L_N(s)-L_N(x))_+^d)  \label{gN}
\end{eqnarray}
and therefore for $1 \le p < \alpha $ sufficiently close to $\alpha $
\begin{eqnarray}
N^{-p H} \E |S^{d,\lambda_N}_{N}(t) - S^{d,\lambda_N}_{N}(s)|^p
&\le&C\|\widetilde g_N(t;\cdot) - \widetilde g_N(s;\cdot)\|_p^p \nn \\
&\le&C\int_0^\infty |(L_N(t)-L_N(s)  + L_N(x))^d  - L_N(x)^d |^p \d x \nn \\
&\le&C\Big((L_N(t)-L_N(s))^{pd}N^{-1} + \int_0^\infty |(L_N(t)-L_N(s)  + x)^d  - x^d |^p \d x \Big) \nn \\
&\le&C\big(L_N(t)-L_N(s)\big)^{1+pd}\label{gN11}
\end{eqnarray}
since $L_N(t) - L_N(s) \ne 0$ implies $L_N(t) - L_N(s)\ge N^{-1}$. This proves
\eqref{tight2}  with $q = 1+ pd  >1$.

Next, let  $ \alpha =2 $ and $ 0 <  H <  \frac{1}{2},$ or  $-1/2< d < 0$.  Then
\eqref{tight2} holds for $S^{d,0}_N$ instead of $S^{d,\lambda_N}_N $ with $q = pH >1 $,
see (\cite{koul}, proof of Prop.~4.4.4).
Hence, it suffices to prove a similar bound for
$R_N(t)$ in \eqref{Sapp}. By Rosenthal's inequality (see the proof \eqref{tight1})
$\E |R_{N}(t) - R_N(s)|^p \le C \E^{p/2} |R_{N}(t) - R_N(s)|^2 $ and hence
$N^{-pH}  \E |R_{N}(t) - R_{N}(s)|^p
\le C \|\widetilde g^0_N(t;\cdot) - \widetilde g^0_N(s;\cdot)\|^p_2 $, where
\begin{eqnarray}\label{g03}
|\widetilde g^0_N(t;x)  - \widetilde g^0_N(s;x)|
&\le&CN^{-d}  \sum_{k=[Ns]+1}^{[Nt]} (k - [Nx])_+^{d-1} \big(1- \e^{-\lambda_N (k - [Nx])_+}\big) \nn \\
&\le&C\begin{cases}(-L_N(x))^{d-1} (L_N(t)-L_N(s)), & x < -1,  \\
C\big((L_N(t)-L_N(x))_+^{d+1}  - (L_N(s)-L_N(x))_+^{d+1}\big)
, &-1 < x < 1,
\end{cases}
\end{eqnarray}
similarly as in  \eqref{gN}.  Hence  $N^{-pH}  \E |R_{N}(t) - R_{N}(s)|^p \le
C \|\widetilde g^0_N(t;\cdot) - \widetilde g^0_N(s;\cdot)\|^p_2 \le
C (L_N(t)-L_N(s))^{1+ p(1+d)} $ follows, proving \eqref{Sapp} and part (ii), too.

\smallskip

\noi (iii) Similarly as in the proof of (ii), let us prove $\|\widetilde g_N(t;\cdot) - g(t;\cdot)\|_p \to 0$, where
\begin{eqnarray} \label{gN3}
\hskip-.5cm &\widetilde g_N(t;x):= \frac{1}{N^d} \sum_{k=1 \vee [Nx]}^{[Nt]}\, b_d (k- [Nx])\, \e^{-\lambda_N (k - [Nx])}, \quad
g(t;x):= \frac{1}{\Gamma(1+d)} h_{H,\alpha,\lambda_*}(t;x),
\end{eqnarray}
see the definition of $h_{H,\alpha,\lambda}$ in \eqref{hdef0}. First, let $d>0 $ or $ H> \frac{1}{\alpha}$. Then
using \eqref{bdk} we obtain the point-wise convergence
\begin{eqnarray} \label{gN4}
\widetilde g_N(t;x)&=&\frac{1}{N \Gamma(d)} \sum_{\frac{1}{N} \vee \frac{[Nx]}{N} \le \frac{k}{N} \le \frac{[Nt]}{N} }
\big(\frac{k}{N} - \frac{[Nx]}{N}\big)^{d-1} \big(1 + \epsilon_{N1}(k,x)\big)
\e^{-\lambda_* (\frac{k}{N} - \frac{[Nx]}{N})(1 + \epsilon_{N2})  }\nn \\
&\to&\frac{1}{\Gamma(d)} \int_0^t (y-x)_+^{d-1}
\e^{-\lambda_*(y-x)_+} d y \ = \   \frac{1}{\Gamma(1+d)} h_{H,\alpha,\lambda_*}(t;x), \quad
\forall \
x \ne 0, t,
\end{eqnarray}
see \eqref{hdef1}, where
$$
\epsilon_{N1}(k,x):= \Gamma(d) (k-[Nx])^{1-d}  b_d (k-[Nx]) - 1 \to 0, \qquad
\epsilon_{N2}
:= (N\lambda_N/\lambda_*) - 1   \to 0
$$
as $N \to \infty, k -[Nx] \to \infty $ and
$|\epsilon_{N1}(k,x)| + |\epsilon_{N2}| < C $ is bounded uniformly in $N, k, x $. Therefore \eqref{gN4} holds
by the dominated convergence theorem. We also have from \eqref{gN3} that
$|\widetilde g_N(t;x)| \le C N^{-d} \sum_{k=1\vee [Nx]}^{[Nt]} (k- [Nx])^{d-1} \e^{-(\lambda_*/2) (k-[Nx])}
\le C \int_0^t (s-x)_+^{d-1} \e^{-(\lambda_*/2)(s-x)_+} \d s  = C h_{d+1/\alpha, \alpha, \lambda_*/2} (t;x) =: \bar g(x)
$ is dominated by an integrable function, see \eqref{hdef1}, with
$\|\bar g \|_p < \infty $. 
This proves \eqref{eq:theorem one 4 part} for $d >0$.

Next, let    $ - \frac{1}{\alpha}< d <0$. Decompose  $\widetilde g_N(t;x)$ in  \eqref{gN3}  as
 $\widetilde g_N(t;x) = \widetilde g_{N1}(x) - \widetilde g_{N2}(t;x)$, where
\begin{eqnarray*} 
\widetilde g_{N1}(x)&:=&N^{-d} \sum_{k=1 \vee [Nx]}^{\infty}\ b_d (k- [Nx]) \, \e^{-\lambda_N (k - [Nx])}, \\
\widetilde g_{N2}(x)&:=&N^{-d} \sum_{[Nt]+1}^{\infty}\ b_d (k- [Nx]) \, \e^{-\lambda_N (k - [Nx])}.
\end{eqnarray*}
Let $ 0< x < t$. First we have
\begin{equation*}
\widetilde g_{N1}(x) 
= \frac{ \sum_{j=0}^{\infty} b_d(j)\, \e^{-\lambda_N j} }{ \lambda_N^{-d} }\ (N \lambda_N)^{-d}\to \lambda_*^{-d}
\end{equation*}
since the last ratio tends to $1$ as $N\to\infty$, see  \eqref{GNlim}.
 We also have
\begin{eqnarray*} 
\widetilde g_{N2}(t;x)&=&\frac{1}{N \Gamma(d)} \sum_{\frac{k}{N} > \frac{[Nt]}{N} }
\big(\frac{k}{N} - \frac{[Nx]}{N}\big)^{d-1} \big(1 + \epsilon_{N1}(k,x)\big)
\e^{-\lambda_* (\frac{k}{N} - \frac{[Nx]}{N})(1 + \epsilon_{N2})  }\nn \\
&\to&\frac{1}{\Gamma(d)} \int_t^\infty (s-x)_+^{d-1} \e^{-\lambda_*(s-x)_+} \d s \ = \ \lambda_*^{-d}
- \frac{1}{\Gamma(1+d)} h_{H,\alpha,\lambda_*}(t;x),
\end{eqnarray*}
see \eqref{hdef2}, similarly to \eqref{gN4}. This proves the point-wise convergence $\widetilde g_N(t;x) \to g(t;x)=\Gamma(1+d)^{-1} h_{H,\alpha,\lambda}(t;x) $ for
$0< x < t$ and the proof for $x < 0$ is similar.  Then
$\|\widetilde g_N(t;\cdot) - g(t;\cdot)\|_p \to 0$ or
\eqref{eq:theorem one 4 part} for $ - \frac{1}{\alpha}< d <0$
follows similarly as in the case $d >0$ above.

Consider the proof of  tightness in (iii). We use the same criterion \eqref{tight2} as
in part (ii). Let first $d >0$. Then for $|x |\le 1$ the bound in
\eqref{gN} and hence $\int_{-1}^1 |\widetilde g_N(t;x) - \widetilde  g_N(s;x)|^p \d x
\le C (L_N(t)-L_N(s))^{1+ pd} $ follows as in \eqref{gN11}. On the other hand,
for $ x<  -1 $ we have
\begin{eqnarray*}
|\widetilde g_N(t;x) - \widetilde g_N(s;x)|
&\le&C \e^{-(\lambda_*/2)|x|}\int_{[Ns/N]}^{[Nt]/N}   \d y \
=\ C\e^{-(\lambda_*/2)|x|}(L_N(t)-L_N(s))  
\end{eqnarray*}
implying $\int_{-\infty}^{-1} |\widetilde g_N(t;x) - \widetilde  g_N(s;x)|^p \d x \le C |L_N(t)-L_N(s)|^p $.
Consequently, \eqref{tight2} for $d>0$ holds with $q = (1+ pd) \wedge p > 1. $
Finally,  \eqref{tight2} for $\alpha = 2,  -1/2 < d < 0$ follows as in case (ii)
since \eqref{g03} holds in the case $\lambda_N = O(1/N) $ as well.
This ends the proof of Theorem~\ref{thm:theorem one}. \hfill $\Box$

\medskip

\noi {\it Proof of Proposition~\ref{propX2}.} (i) By stationarity, $N^{-1}\E \sum_{t=1}^N X^2_{d,\lambda_N}(t) = \E X^2_{d,\lambda_N}(0).  $
Let us first prove the convergence of expectations:
\begin{eqnarray}
\E X^2_{d,\lambda_N}(0)
&\to&\frac{\Gamma(1-2d)}{\Gamma^2(1-d)}, \hskip1cm  d < 1/2,  \label{Esum1}  \\
\lambda_N^{2d-1} \E X^2_{d,\lambda_N}(0)
&\to&\frac{\Gamma(d-1/2)}{2\sqrt{\pi}\,\Gamma(d)},  \qquad d > 1/2, \label{Esum2} \\
|\log \lambda_N|^{-1}
\E X^2_{d,\lambda_N}(0)&\to&\frac{1}{\pi}, \hskip2cm d=1/2, \label{Esum3}
\end{eqnarray}
as $N \to \infty$. Since $\E X^2_{d,\lambda_N}(0) = \sum_{k=0}^\infty \e^{-2\lambda_N k} \omega^2_{-d}(k), $
with $\omega_{-d}(k) \equiv \omega_{-d}(k)$ defined in \eqref{omega0},
and
\begin{equation*}
\sum_{k=0}^n \omega^2_{-d}(k) \ \sim \
\begin{cases} (1/(2d-1)\Gamma^2(d)) n^{2d-1}, &d> 1/2, \\
 (1/\Gamma^2(1/2)) \log (n),  &d=1/2, \\
\sum_{k=0}^\infty \omega^2_{-d}(k) = \Gamma(1-2d)/\Gamma^2(1-d), &d < 1/2
\end{cases}
\end{equation*}
as  $n \to \infty$, the convergences in \eqref{Esum1}-\eqref{Esum3} follows from the Tauberian theorem in \cite{fell1966} \
used in the proof of Theorem~4.3 (i) above.

With \eqref{Esum1}-\eqref{Esum3} in mind,
\eqref{sum1}-\eqref{sum3} follow from
\begin{eqnarray}\label{QN}
Q_N \equiv N^{-1}\sum_{t=1}^N \Big[X^2_{d,\lambda_N}(t) - \E X^2_{d,\lambda_N}(t)  \big]
&=&\begin{cases}o_p(1), &d < 1/2, \\
o_p(\lambda_N^{1-2d}), &d>1/2, \\
o_p(|\log \lambda_N|), &d=1/2.
\end{cases}
\end{eqnarray}
Let $\omega_{-d,\lambda}(k) := \omega_{-d}(k) \e^{-\lambda k} $.
We have $Q_N = Q_{N1} + Q_{N2} $, where
\begin{eqnarray*}
Q_{N1}&=&N^{-1} \sum_{s \le N} (\zeta^2(s) - \E \zeta^2(s)) \sum_{t=1\vee s}^N \omega^2_{-d, \lambda_N}(t-s), \\
Q_{N2}&=&N^{-1}\sum_{s_2 < s_1 \le N} \zeta (s_1)  \sum_{t=1 \vee s_1}^N
 \omega_{-d,\lambda_N}(t-s_1) \omega_{-d,\lambda_N}(t-s_2) \zeta(s_2).
\end{eqnarray*}
Note $Q_{Ni}, i=1,2 $ are sums of martingale differences. We shall use the well-known moment inequality for sums of martingale differences:
\begin{equation} \label{sum5}
\E |\sum_{i\ge 1} \xi_i |^\alpha \le 2 \sum_{i\ge 1} \E |\xi_i|^\alpha
\end{equation}
see e.g. (\cite{koul}, Prop.~2.5.2),
which is  valid for any
$1\le \alpha \le 2$ and any
sequence $\{\xi_i, i \ge 1\}$ with $\E |\xi_i|^\alpha < \infty,
\E [\xi_i |\xi_j, 1\le j < i ] = 0, i \ge 1$.  \\
First, let $d > 1/2 $. Using $\E |\zeta(0)|^p < \infty, 2< p < 4 $ and
\eqref{sum5} with $\alpha = p/2 $
 we obtain
\begin{eqnarray*}
\E |Q_{N1}|^{p/2}
&\le&CN^{-p/2}\sum_{s \le N} \Big|\sum_{t=1\vee s}^N \omega^2_{-d,\lambda_N}(t-s) \Big|^{p/2}\nn \\
&\le&CN^{-p/2} \sum_{s \le N} \Big(\sum_{t=1\vee s}^N (t-s)^{2(d-1)}_+ \e^{-2\lambda_N(t-s)} \Big)^{p/2} \nn \\
&\le&CN^{-p/2}\Big\{\int_N^\infty \d s \Big(\int_0^N (t+s)^{2(d-1)} \e^{-2\lambda_N(t+s)} \d t \Big)^{p/2}
+ N \Big(\int_0^N t^{2(d-1)} \e^{-2\lambda_N t} \d t \Big)^{p/2} \Big\} \\
&\equiv&CN^{-p/2} \big\{I_{N1} + I_{N2}\big\}.
\end{eqnarray*}
Here, $I_{N2} \le N \Big(\int_0^\infty t^{2(d-1)} \e^{-2\lambda_N t} \d t \Big)^{p/2} = C N \lambda_N^{(p/2)(1-2d)} $ and
$I_{N1} \le  C N^{p/2} \int_N^\infty s^{(d-1)p} \e^{- \lambda_N p s} \d  s  \le C N^{p/2} \lambda_N^{(1-d)p -1}. $ Therefore,
\begin{equation}
\E |Q_{N1}|^{p/2} \le C\big(\lambda_N^{(1-d)p -1}  +
N^{1-p/2} \lambda_N^{(p/2)(1-2d)}\big).  
\label{QN1}
\end{equation}
Next,
\begin{eqnarray}
\E |Q_{N2}|^{2}
&\le&CN^{-2}\sum_{s_2 < s_1 \le N} \E \big(\sum_{t=1 \vee s_1}^N  \omega_{-d,\lambda_N}(t-s_1)
\omega_{-d,\lambda_N}(t-s_2) \zeta(s_2)\big)^2 \nn \\
&\le&CN^{-2}\sum_{s_2 < s_1 \le N} \sum_{t=1 \vee s_1}^N  \omega^2_{-d,\lambda_N}(t-s_1)
\omega^2_{-d,\lambda_N}(t-s_2) \nn \\
&\le&CN^{-2}\sum_{t=1}^N \big(\sum_{s \le t}  \omega^2_{-d,\lambda_N}(t-s) \big)^2 \nn \\
&\le&CN^{-2} \int_0^N \d t \Big( \int_{-\infty}^t  (t-s)^{2(d-1)} \e^{-2\lambda_N(t-s)} \d s \Big)^2 \nn \\
&=&CN^{-1}\Big( \int_0^{\infty} s^{2(d-1)} \e^{-2\lambda_Ns} \d s \Big)^2 \
\le\ CN^{-1} \lambda_N^{2(1-2d)}. \label{QN2}
\end{eqnarray}
Since $(1-d)p-1 >  (p/2)(1-2d)$ and $1-(p/2) < 0$,
\eqref{QN1}  and \eqref{QN2}  prove \eqref{QN} for $d >  1/2$. \\
Next, let $d < 1/2$. Then similarly as above we obtain
\begin{eqnarray}
\E |Q_{N1}|^{p/2}
&\le&CN^{-p/2} \sum_{s \le N} \Big(\sum_{t=1\vee s}^\infty (t-s)^{2(d-1)}_+  \Big)^{p/2} \nn \\
&=&CN^{-p/2}\Big\{N \big(\sum_{t=1}^N t^{2(d-1)} \big)^{p/2} + \sum_{s \ge N} \big(\sum_{t=1}^N (s+t)^{2(d-1)} \big)^{p/2} \Big\} \nn \\
&\le&CN^{1-p/2}  + C N^{-p/2} \sum_{s \ge N}  (N s^{2(d-1)})^{p/2}\ \le\  CN^{1-p/2} \label{QN3}
\end{eqnarray}
and
\begin{eqnarray}
\E |Q_{N2}|^{2}
&\le&CN^{-2}\sum_{s_2 < s_1 \le N} \sum_{t=1 \vee s_1}^N  \omega^2_{-d,\lambda_N}(t-s_1) \omega^2_{-d,\lambda_N}(t-s_2) \nn \\
&\le&CN^{-2} \sum_{t=1}^N \Big( \sum_{s < t}  (t-s)^{2(d-1)} \Big)^2 \ \le\ CN^{-1}. \label{QN4}
\end{eqnarray}
\eqref{QN3}  and \eqref{QN4}  prove \eqref{QN} for $d <  1/2$. \\
Finally, let $d=1/2$. Then since $p>2 $
\begin{eqnarray}
\E |Q_{N1}|^{p/2}
&\le&CN^{-p/2} \sum_{s \le N} \Big(\sum_{t=1\vee s}^\infty (t-s)^{-1}_+  \Big)^{p/2} \nn \\
&=&CN^{-p/2}\Big\{N \big(\sum_{t=1}^N t^{-1} \big)^{p/2} + \sum_{s \ge N} \big(\sum_{t=1}^N (s+t)^{-1} \big)^{p/2} \Big\} \nn \\
&\le&CN^{1-p/2} (\log N)^{p/2}  + C N^{-p/2} \sum_{s \ge N}  (N s^{-1})^{p/2} \nn \\
&\le&CN^{1-p/2} (\log N)^{p/2} \ = \ o(1) \label{QN5}
\end{eqnarray}
while
\begin{eqnarray}
\E |Q_{N2}|^{2}
&\le&CN^{-2} \sum_{t=1}^N \Big( \sum_{s < t}  (t-s)^{-1} \e^{-\lambda_N (t-s)}\Big)^2 \
\le\ CN^{-1}(\log \lambda_N)^2. \label{QN6}
\end{eqnarray}
\eqref{QN5}  and \eqref{QN6}  prove \eqref{QN} for $d =  1/2$. Proposition~\ref{propX2} is proved.
\hfill $\Box$


\bigskip


\footnotesize

\begin{thebibliography}{99}




\bibitem{abramowitz}
Abramowitz, M., Stegun, I.: Handbook of mathematical functions, ninth edition, Dover, New York (1965).







\bibitem{astrauskas} Astrauskas, A. (1983).
Limit theorems for sums of linearly generated random variables.
{\it Lithuanian J. Math.} {\bf 23}\  127--134.






\bibitem{balan} Balan, R., Jakubowski, A. and Louhichi, S.  (2016).
Functional convergence of linear processes with heavy-tailed innovations.
{\it J. Theoret. Probab.} {\bf 29}\  491--526.











\bibitem{Bill}
Billingsley, P. (1968). {\it Convergence of Probability Measures}. New York: Wiley.






\bibitem{BrockwellDavisTSTM}
Brockwell, P.J. and Davis, R.A. (1991).
{\it Time Series: Theory and Methods}, 2nd ed.. New York: Springer.








\bibitem{Dickey1}
Dickey, D. and Fuller, W. (1979). Distribution of the estimators
for autoregressive time series with a unit root. {\it JASA}
{\bf 74}\  427--431.





\bibitem{fell1966} Feller, W. (1966). {\it An Introduction to Probability Theory and Its Applications, vol. 2.}
New York: Wiley.




\bibitem{giraitis}
Giraitis, L., Kokoszka, P. and Leipus, R. (2000). Stationary ARCH models: dependence structure
and central limit theorem. {\it Econometric Theory} {\bf 16}\  3--22.


\bibitem{giraitis2003a}
Giraitis, L., Kokoszka, P., Leipus, R. and Teyssi\`ere, G. (2003).
Rescaled variance and related tests
for long memory in volatility and levels. {\it J. Econometrics} {\bf 112}\  265--294.





\bibitem{koul}
Giraitis, L., Koul. H.L. and Surgailis, D. (2012). {\it Large Sample Inference for Long Memory Processes.}
London: Imperial College Press.

\bibitem{Gradshteyn}
Gradshteyn, I.S. and Ryzhik, I.M. (2000).
{\it Tables of Integrals and Products.} 6th edition. New York: Academic Press.



\bibitem{Granger}
Granger, C.W.J. and Joyeux, R. (1980).
An introduction to long-memory time series models and fractional differencing.
{\it J. Time Series Anal.}  {\bf 1}\  15--29.




\bibitem{Hosking}
Hosking, J.R.M. (1981). Fractional differencing. {\it Biometrika} {\bf 68} \ 165--176.




\bibitem{kasahara} Kasahara, Y. and Maejima, M. (1988). Weighted sums of i.i.d. random variables
attracted to integrals of stable processes. {\it Probab. Theory Relat. Fields } {\bf 78}\  75--96.



\bibitem{kokoszka} Kokoszka, P.S. and Taqqu, M.S. (1995). Fractional ARIMA with stable innovations.
{\it Stochastic Process.  Appl.} {\bf 60} \ 19--47.












\bibitem{Lavancier} Lavancier, F., Leipus, R.,  Philippe, A. and Surgailis, D. (2013). Detection of
non-constant long memory parameter.
{\it Econometric Theory} {\bf 29}\  1009--1056.



\bibitem{Lo}
Lo, A. (1991). Long-term memory in stock market prices. {\it Econometrica} {\bf 59}\ 1279--1313.
















\bibitem{Meerschaertsabzikar}
Meerschaert, M.M. and Sabzikar, F. (2013).
Tempered fractional Brownian motion. {\it  Statist. Probab. Lett.} {\bf 83} \ 2269--2275.



\bibitem{Meerschaertsabzikar2} Meerschaert, M.M. and Sabzikar, F. (2016). Tempered fractional stable motion.
{\it J. Theoret. Probab.} {\bf 29}\ 681--706.


\bibitem{Meerschaertsabzikarkumarzeleki}
Meerschaert, M.M., Sabzikar, F., Phanikumar, M.S. and Zeleke, A. (2014).
Tempered fractional time series model for turbulence in geophysical flows.
{\it J. Stat. Mech. Theory  Exp.} {\bf 2014} \  P09023.











\bibitem{Phillips}
Phillips, P.C.B. (1987). Time series regression with a unit root. {\it Econometrica} {\bf 55}\ 277--301.











\bibitem{SurgailisFarzadTFSMII}
Sabzikar, F. and Surgailis, D. (2017). Tempered fractional Brownian and stable motions of second kind.
Preprint. Available on http://arxiv.org/abs/1702.07258





\bibitem{SamorodnitskyTaqqu}
Samorodnitsky, G. and Taqqu, M.S.  (1996). {\it Stable Non-Gaussian Random Processes: Stochastic Models with Infinite Variance.}
Boca Raton etc: Chapman and Hall.




\bibitem{sowell2}
Sowell, F. (1990). The fractional unit root distribution. {\it Econometrica} {\bf 58}\   495--505.






















\end{thebibliography}
\end{document}